\documentclass[a4,12pt]{article}
\usepackage{amsmath}
\usepackage{amssymb}
\usepackage{amsthm}
\newtheorem{prop}{Proposition}
\theoremstyle{definition}
\newtheorem*{remarks}{Remarks}
\newtheorem*{note}{Note}
\newtheorem*{reference}{References for Remarks}
\newtheorem*{corrections}{Corrections to [I1] and [I2]}

\begin{document}
\title{On generalizations of the sum formula for multiple zeta values}
\author{Masahiro Igarashi}
\date{}
\maketitle
\begin{abstract}
We prove some generalizations of the sum formula for multiple zeta values 
by using Hiroyuki Ochiai's method of proving the sum formula. 
\end{abstract}
\begin{note} 
The present paper is an English summary of my master's thesis written in Japanese: Masahiro Igarashi, ``On generalizations of the sum formula for multiple zeta values", 
Master's thesis, Graduate School of Mathematics, Nagoya University, Japan, 
February 3, 2007 (in Japanese). In this summary, 
I corrected an error. For details, see page 7. 
In addition, I added some related remarks at the end. 
\end{note}
\section{Introduction}
The multiple zeta value (MZV for short) is defined by the multiple series
\begin{equation*}
 \zeta(k_1,\ldots,k_n) := \sum_{0< m_1< \cdots <m_n<\infty}\frac{1}{m_1^{k_1} \cdots m_n^{k_n}},
\end{equation*}
where $k_1,\ldots ,k_n\in\mathbb{N}:=\{1,2,3,\ldots\}$ with $k_n\ge2$ (Euler \cite{eu}, Hoffman \cite{h1} and Zagier \cite{z}; 
see also the surveys \cite{ak}, \cite{kan} and \cite{k2}). 
It is known that MZVs satisfy various relations over $\mathbb{Q}$. 
The sum formula is one of the well-known $\mathbb{Q}$-linear 
relations among MZVs. 
It was first proved independently by A. Granville \cite{gr} and 
D. Zagier (unpublished). For the double and triple zeta values 
and the sum conjecture, see Euler \cite{eu}; Hoffman and Moen \cite{hm}; 
Hoffman \cite{h1}, Markett \cite{m}, respectively. 
For alternative proofs and generalizations of the sum formula, 
see, e.g., \cite{ak}, \cite{ho}, \cite{kan}, \cite{o1}, \cite{o2}, \cite{oz}. 
In the present paper, we prove the following three generalizations of the sum formula: 
\begin{prop}
We have
\begin{equation*}
\zeta(k;\alpha)
=
\sum_{\begin{subarray}{c}k_1+\cdots+k_{n}=k\\
       k_i\in\mathbb{N}, k_n\ge2
      \end{subarray}}
\sum_{0\le m_1< \cdots <m_n<\infty}\frac{(\alpha)_{m_1}}{{m_1}!}\frac{{m_n}!}{(\alpha)_{m_n}}
\frac{1}{(m_1+\alpha)^{k_1} \cdots (m_n+\alpha)^{k_n}}
\end{equation*}
for all $k, n\in\mathbb{N}$ such that $n<k$ and $\alpha\in\mathbb{C}$ with 
$\mathrm{Re}(\alpha)>0$, 
where $\zeta(s;\alpha):=\sum_{m=0}^{\infty}(m+\alpha)^{-s}$, the Hurwitz zeta function, and $(a)_m$ denotes the Pochhammer symbol, i.e., 
$(a)_m=a(a+1)\cdots(a+m-1)$ $(m\in\mathbb{N})$ and 
$(a)_0=1$. 
\end{prop}
For properties of the Hurwitz zeta function, see \cite{sc} for instance. 
\begin{prop}
We have
\begin{equation}
\begin{aligned}
&\sum_{\begin{subarray}{c}k_1+\cdots+k_n=k\\
k_i\in\mathbb{N},k_n\ge2
\end{subarray}}
\zeta(k_1,\ldots,k_n;\alpha)\\
=
&\frac{1}{(k-n-1)!}\sum_{l=0}^{\infty}\frac{1}{(l+1)^n}
\frac{\partial^{k-n-1}}{{\partial}X^{k-n-1}}
\Biggl(
\frac{(1-X)_l}{(\alpha-X)_{l+1}}
\Biggr)
{\Biggl|}_{X=0}
\end{aligned}
\end{equation}
for all $k, n\in\mathbb{N}$ such that $n<k$ and $\alpha\in\mathbb{C}$ with 
$\mathrm{Re}(\alpha)>0$, where 
\begin{equation*}
\zeta(s_1,\ldots,s_n;\alpha)
:=
\sum_{0{\le}m_1<\cdots<m_n<\infty}
\frac{1}{(m_1+\alpha)^{s_1}\cdots(m_n+\alpha)^{s_n}},
\end{equation*}
the multiple Hurwitz zeta function.
\end{prop}
For Proposition 1, we will prove the following extension: 
\begin{prop}
We have
\begin{equation}
\begin{aligned}
&\sum_{l=0}^{\infty}\frac{1}{(l+\alpha)^n(l+\beta)^m}\\
=
&\sum_{\begin{subarray}{c}k_1+\cdots+k_{n}=m+n-1\\
       k_i\in\mathbb{N}
      \end{subarray}}
\sum_{0\le m_1< \cdots <m_n<\infty}\frac{(\alpha)_{m_1}}{{m_1}!}\frac{{m_n}!}{(\alpha)_{m_{n}+1}}
\frac{1}{(m_1+\beta)^{k_1} \cdots (m_n+\beta)^{k_n}}
\end{aligned}
\end{equation}
for all $m, n\in\mathbb{N}$ and $\alpha,\beta\in\mathbb{C}$ with 
$\mathrm{Re}(\alpha),\mathrm{Re}(\beta)>0$. 
\end{prop}
The case $\alpha=\beta$ is Proposition 1. All the above propositions are 
generalizations of the sum formula for MZVs: for $\alpha=\beta=1$, 
they each become it. To prove Propositions 2 and 3, we use Hiroyuki Ochiai's method of proving the sum formula. Though Ochiai's proof is unpublished, 
it can be found in \cite[pp.~17--20]{ak} and \cite[pp.~60--61]{kan}.
\section{Proofs of Propositions 2 and 3}
\begin{proof}[Proof of Proposition $3$]
Let $\alpha,\beta\in\mathbb{C}$ with $\mathrm{Re}(\alpha), \mathrm{Re}(\beta)>0$ and $n\in\mathbb{N}$. 
Using the well-known formula
\begin{equation*}
\int_{0}^{1}t^{a-1}(1-t)^{b-1}\,\mathrm{d}t =\frac{\Gamma(a)\Gamma(b)}{\Gamma(a+b)}
\end{equation*}
($a, b \in\mathbb{C}$ with $\mathrm{Re}(a), \mathrm{Re}(b)>0$), 
we have the following iterated integral representation of a generating function of the left-hand side of (2):
\begin{equation}
\begin{aligned}
&\sum_{m=0}^{\infty}\left(\sum_{l=0}^{\infty}\frac{1}{(l+\alpha)^n(l+\beta)^{m+1}}\right)X^m\\
=
&\sum_{l=0}^{\infty}\frac{1}{(l+\alpha)^n(l+\beta-X)}\\
=
&\int_{0}^{1}\frac{(1-t_0)^{-X+\beta-1}}{t_0^\alpha}\,\mathrm{d}t_0 
\int_{0}^{t_0}\frac{\mathrm{d}t_1}{t_1} \cdots 
\int_{0}^{t_{n-2}}\frac{\mathrm{d}t_{n-1}}{t_{n-1}}
\int_{0}^{t_{n-1}}\frac{t_n^{\alpha-1}}{(1-t_n)^{-X+\beta}}\,\mathrm{d}t_n\\
=
&\idotsint\displaylimits_{0<t_n<\cdots< t_0<1}
\frac{(1-t_0)^{-X+\beta-1}t_n^{\alpha-1}\,\mathrm{d}t_0 \cdots \mathrm{d}t_n}
{t_0^{\alpha}{t_1}\cdots{t_{n-1}}(1-t_n)^{-X+\beta}}
\end{aligned}
\end{equation}
for $X\in\mathbb{C}$ such that $|X|<\mathrm{Re}(\beta)$. 
Making the change of variables
\begin{equation*}
t_i \,\mapsto\, 1-t_{n-i} \quad (i=0,1,\ldots,n)
\end{equation*}
(see Zagier \cite[p.~510]{z}) to the integral, we have 
\begin{equation}
\begin{aligned}
&\idotsint\displaylimits_{0<t_n<\cdots< t_0<1}
\frac{(1-t_0)^{-X+\beta-1}t_n^{\alpha-1}\,\mathrm{d}t_0 \cdots \mathrm{d}t_n}
{t_0^{\alpha}{t_1}\cdots{t_{n-1}}(1-t_n)^{-X+\beta}}\\
=
&\idotsint\displaylimits_{0<t_n<\cdots< t_0<1}
\frac{t_n^{-X+\beta-1}(1-t_0)^{\alpha-1}\,\mathrm{d}t_0 \cdots \mathrm{d}t_n}
{(1-t_n)^{\alpha}{(1-t_{n-1})}\cdots{(1-t_1)}t_0^{-X+\beta}}
\end{aligned}
\end{equation}
for $|X|<\mathrm{Re}(\beta)$. 
Further, the right-hand side of (4) can be calculated as follows: 
\begin{equation}
\begin{aligned}
&\idotsint\displaylimits_{0<t_n<\cdots <t_0<1}
\frac{t_n^{-X+\beta-1}(1-t_0)^{\alpha-1}\,\mathrm{d}t_0 \cdots \mathrm{d}t_n}
{(1-t_n)^{\alpha}{(1-t_{n-1})}\cdots{(1-t_1)}t_0^{-X+\beta}}\\
=&\int_{0}^{1}\frac{(1-t_0)^{\alpha-1}}{t_0^{-X+\beta}}\,\mathrm{d}t_0 
\int_{0}^{t_0}\frac{\mathrm{d}t_1}{1-t_1} \cdots 
\int_{0}^{t_{n-2}}\frac{\mathrm{d}t_{n-1}}{1-t_{n-1}}
\int_{0}^{t_{n-1}}\frac{t_n^{-X+\beta-1}}{(1-t_n)^{\alpha}}\,\mathrm{d}t_n\\
=&\int_{0}^{1}\frac{(1-t_0)^{\alpha-1}}{t_0^{-X+\beta}}\,\mathrm{d}t_0 
\int_{0}^{t_0}\frac{\mathrm{d}t_1}{1-t_1} \cdots 
\int_{0}^{t_{n-2}}\frac{\mathrm{d}t_{n-1}}{1-t_{n-1}}
\sum_{l_1=0}^{\infty}\frac{(\alpha)_{l_1}}{{l_1}!}\frac{t_{n-1}^{l_1+\beta-X}}{l_1+\beta-X}\\
=&\cdots\cdots\\
=&\sum_{l_1 ,\ldots ,\l_n\ge0}
\frac{(\alpha)_{l_1}}{l_1!}\frac{1}{(l_1+\beta-X)\cdots(l_1+\cdots+l_n+n-1+\beta-X)}\\
&\hspace*{1cm}\times\int_{0}^{1}(1-t_0)^{\alpha-1}t_{0}^{l_1+\cdots+l_n+n-1}\,\mathrm{d}t_0\\
=
&\sum_{l_1 ,\ldots ,\l_n\ge0}
\frac{(\alpha)_{l_1}}{l_1!}\frac{1}{(l_1+\beta-X)\cdots(l_1+\cdots+l_n+n-1+\beta-X)}\\
&\hspace*{1cm}\times\frac{\Gamma(\alpha)\Gamma(l_1+\cdots+l_n+n)}{\Gamma(\alpha+l_1+\cdots+l_n+n)}\\
=&\sum_{0\le m_1< \cdots <m_n<\infty}\frac{(\alpha)_{m_1}}{{m_1}!}\frac{{m_n}!}{(\alpha)_{m_{n}+1}}
\frac{1}{(m_1+\beta-X) \cdots (m_n+\beta-X)}\\
=&\sum_{m=0}^{\infty}\left(\sum_{\begin{subarray}{c}k_1+\cdots+k_{n}=m+n\\
       k_i\in\mathbb{N}
      \end{subarray}}
\sum_{0\le m_1< \cdots <m_n<\infty}\frac{(\alpha)_{m_1}}{{m_1}!}\frac{{m_n}!}{(\alpha)_{m_{n}+1}}
\left\{\prod_{j=1}^{n}\frac{1}{(m_j+\beta)^{k_j}}\right\}\right){X^m}
\end{aligned}
\end{equation}
for $|X|<\mathrm{Re}(\beta)$. 
Hence, combining all the identities (3), (4) and (5), we obtain (2). 
\end{proof}
\begin{proof}[Proof of Proposition $2$] 
The proof is similar to that of Proposition 3. Let $\alpha\in\mathbb{C}$ with $\mathrm{Re}(\alpha)>0$ and $n\in\mathbb{N}$. 
By using the same calculus as in the proof of (3), we also have the following iterated integral 
representation of a generating function of the left-hand side of (1):
\begin{equation}
\begin{aligned}
&\sum_{m=0}^{\infty}\left(\sum_{\begin{subarray}{c}k_1+\cdots+k_n=m+n+1\\
k_i\in\mathbb{N},k_n\ge2
\end{subarray}}
\zeta(k_1,\ldots,k_n;\alpha)\right)X^m\\
=
&\idotsint\displaylimits_{0<t_n<\cdots< t_0<1}
\frac{t_0^{X-1}t_n^{\alpha-X-1}\,\mathrm{d}t_0 \cdots \mathrm{d}t_n}
{(1-t_1)\cdots(1-t_n)}
\end{aligned}
\end{equation}
for $X\in\mathbb{C}$ such that $|X|<\mathrm{Re}(\alpha)$. 
Similarly, making the change of variables
$t_i \,\mapsto\, 1-t_{n-i}$ ($i=0,1,\ldots,n$) to the integral gives 
\begin{equation}
\begin{aligned}
&\idotsint\displaylimits_{0<t_n<\cdots< t_0<1}
\frac{t_0^{X-1}t_n^{\alpha-X-1}\,\mathrm{d}t_0 \cdots \mathrm{d}t_n}
{(1-t_1)\cdots(1-t_n)}\\
=
&\idotsint\displaylimits_{0<t_n<\cdots< t_0<1}
\frac{(1-t_0)^{\alpha-X-1}\,\mathrm{d}t_0 \cdots \mathrm{d}t_n}
{(1-t_n)^{1-X}t_{n-1}\cdots{t_0}}
\end{aligned}
\end{equation}
for $|X|<\mathrm{Re}(\alpha)$. 
Further, the same calculus as in the proof of (5) gives the following series 
expression of the right-hand side of (7): 
\begin{equation}
\begin{aligned}
&\idotsint\displaylimits_{0<t_n<\cdots< t_0<1}
\frac{(1-t_0)^{\alpha-X-1}\,\mathrm{d}t_0 \cdots \mathrm{d}t_n}
{(1-t_n)^{1-X}t_{n-1}\cdots{t_0}}\\
=
&\sum_{l=0}^{\infty}\frac{(1-X)_l}{(\alpha-X)_{l+1}}\frac{1}{(l+1)^n}
\end{aligned}
\end{equation}
for $|X|<\mathrm{Re}(\alpha)$. 
We see that the factor ${(1-X)_l}/{(\alpha-X)_{l+1}}$ has the expansion 
\begin{equation*}
\frac{(1-X)_l}{(\alpha-X)_{l+1}}
=
\sum_{m=0}^{\infty}\frac{1}{m!}
\frac{\partial^m}{{\partial}X^m}
\Biggl(
\frac{(1-X)_l}{(\alpha-X)_{l+1}}
\Biggr)
{\Biggl|}_{X=0}
X^m
\quad \mathrm{for}\,\, |X|<\mathrm{Re}(\alpha).
\end{equation*}
Hence, substituting this into the right-hand side of (8), we have 
\begin{equation}
\begin{aligned}
&\sum_{l=0}^{\infty}\frac{(1-X)_l}{(\alpha-X)_{l+1}}\frac{1}{(l+1)^n}\\
=
&\sum_{l=0}^{\infty}\sum_{m=0}^{\infty}
\frac{1}{m!}\frac{\partial^m}{{\partial}X^m}
\Biggl(
\frac{(1-X)_l}{(\alpha-X)_{l+1}}
\Biggr)
{\Biggl|}_{X=0}
\frac{X^m}{(l+1)^n}
\end{aligned}
\end{equation}
for $|X|<\mathrm{Re}(\alpha)$. 
Here we note that the double series on the right-hand side of (9) converges absolutely for $X\in\mathbb{C}$ such that $|X|<r$, 
where $r$ is a fixed real number such that $0<r<\mathrm{Re}(\alpha)/2$. 
(In my master's thesis, I imposed not $|X|<r$ but $|X|<\mathrm{Re}(\alpha)/2$ in order to prove the absolute convergence of the double series. 
I think that the latter condition is incorrect.) 
Thus we have 
\begin{equation}
\begin{aligned}
&\idotsint\displaylimits_{0<t_n<\cdots< t_0<1}
\frac{(1-t_0)^{\alpha-X-1}\,\mathrm{d}t_0 \cdots \mathrm{d}t_n}
{(1-t_n)^{1-X}t_{n-1}\cdots{t_0}}\\
=
&\sum_{m=0}^{\infty}
\left(\sum_{l=0}^{\infty}
\frac{1}{(l+1)^n}
\frac{\partial^m}{{\partial}X^m}
\Biggl(
\frac{(1-X)_l}{(\alpha-X)_{l+1}}
\Biggr)
{\Biggl|}_{X=0}
\right)\frac{X^m}{m!}
\end{aligned}
\end{equation}
for $|X|<r$. 
Combining all the identities (6), (7) and (10), we obtain (1). 
\end{proof}

\pagebreak
\begin{remarks}[\today] 
I give some remarks related to my master's thesis. 
\par  
(i) My master's thesis and my subsequent paper \cite{I2} got benefits from Ochiai's proof of the sum formula for MZVs. In particular, I discovered the one- and two-parameter multiple series in Propositions 1 and 3 and their iterated integral representations by studying Ochiai's proof. Ochiai's iterated integral representation of a generating function played 
an essential role for these discoveries. Ochiai's proof can be found also in \cite[pp.~19--22]{AK}. 
D. Zagier's unpublished proof of the sum formula can also be found in \cite[pp.~18--19]{AK}. 
Zagier also used an iterated integral representation of a generating function. 
One of the differences between Ochiai's iterated integral representation and 
Zagier's is that Ochiai's has an exponential parameter, whereas Zagier's a polynomial parameter. 
From this difference, I conclude that it is impossible to discover my one- and two-parameter multiple series and their iterated integral representations by examining Zagier's proof and its generalization of \cite{O}. 
Here Ochiai's proof has an advantage. My multiple series and their iterated integral representations are highly non-trivial. I proved Propositions 1--3 in 2006. In 2006, I talked on Proposition 1 and its proof at my then academic adviser's seminar. On February 3, 2007, my master's thesis passed the final examination. 
From some facts, I concluded in February 2007 that a certain person had already been aware of the contents of my master's thesis at the workshop ``Zeta Wakate Kenky\={u}sh\={u}kai" held at Graduate School of Mathematics, Nagoya University, Japan, February 17--18, 2007. 
\par
(ii) Multiple series similar to those in Propositions 1 and 3 
appear in the study of the $\varepsilon$-expansion of hypergeometric functions; 
see, e.g., \cite{DK}, \cite{K}, \cite{KK}, \cite{KV}, \cite{KWY}. 
For other similar multiple series, see \cite{I2} and its preprint arXiv:0908.2536v7. 
\par 
(iii) In \cite[pp.~752--753]{I3}, I gave an explicit expression of the right-hand side of (1) by using a calculus of mine for the Pochhammer symbol $(a)_m$. 
\end{remarks}
\begin{corrections}
I make the following corrections, because the date when my master's thesis 
passed the final examination is February 3, 2007. 
Page 517, line 12 from the bottom of \cite{I1} and page 578, line 23 from the 
bottom of \cite{I2}: ``March 2007" should be ``February 3, 2007". 
\end{corrections}
\begin{flushleft}
Nagoya, Japan\\
\textit{E-mail address}: masahiro.igarashi2018@gmail.com
\end{flushleft}
\end{document}